\documentclass[11pt]{amsart}
\input{epsf}
\usepackage{amssymb,latexsym,float}
\topskip  \numberwithin{equation}{section}
\newtheorem{theorem}{Theorem}[section]
\newtheorem{mtheorem}{Main Theorem}[section]

\newtheorem{remark}[theorem]{Remark}
\newtheorem{lemma}[theorem]{Lemma}
\newtheorem{fact}[theorem]{Fact}

\begin{document}

\pagenumbering{arabic} \pagestyle{headings}
\def\sof{\hfill\rule{2mm}{2mm}}
\def\ls{\leq}
\def\gs{\geq}
\def\SS{\mathcal S}
\def\qq{{\bold q}}
\def\txx{{\frac1{2\sqrt{x}}}}
\def\tx{{\left(\txx\right)}}
\def\ttx{\left({{\frac{1}{2x}}}\right)}
\def\Bn{\mathcal{P}_n}
\def\mn{\mbox{-}}
\def\vn{\varnothing}

\title{ {\sc  Counting peaks at height $k$ in a Dyck path}}

\author{Toufik Mansour}
\maketitle
\begin{center}{LaBRI, Universit\'e Bordeaux 1,
              351 cours de la Lib\'eration\\
              33405 Talence Cedex, France\\[4pt]
        {\tt toufik@labri.fr} }
\end{center}
\markboth{{\normalsize Toufik Mansour}} {{\normalsize Counting
peaks at height $k$ in a Dyck path}}

\section*{Abstract}
A Dyck path is a lattice path in the plane integer lattice
$\mathbb{Z}\times\mathbb{Z}$ consisting of steps $(1,1)$ and
$(1,-1)$, which never passes below the $x$-axis. A peak at height
$k$ on a Dyck path is a point on the path with coordinate $y=k$
that is immediately preceded by a $(1,1)$ step and immediately
followed by a $(1,-1)$ step. In this paper we find an explicit
expression for the generating function for the number of Dyck
paths starting at $(0,0)$ and ending at $(2n,0)$ with exactly $r$
peaks at height $k$. This allows us to express this function via
Chebyshev polynomials of the second kind and the generating
function for the Catalan numbers.

\noindent{\bf Keywords}: Dyck paths, Catalan numbers, Chebyshev
polynomials.
\section{Introduction and main results}
The {\em Catalan sequence} is the sequence
$$\{C_n\}_{n\geq0}=\{1,1,2,5,14,132,429,1430,\dots\},$$
where $C_n=\frac{1}{n+1}\binom{2n}{n}$ is called the $n$th {\em
Catalan number}. The generating function for the Catalan numbers
is denoted by $C(x)=\frac{1-\sqrt{1-4x}}{2x}$. The Catalan numbers
provide a complete answer to the problem of counting certain
properties of more than $66$ different combinatorial structures
(see Stanley \cite[Page~219 and Exercise~6.19]{St}). The structure
of use to us in the present paper is Dyck paths.

\begin{center}
\begin{figure}[H]
\hspace*{1.25truecm}
\epsfxsize=360.0pt 
\epsffile{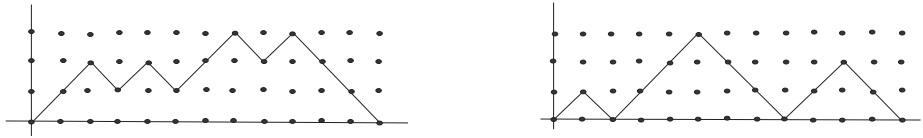} \caption{Two Dyck paths.} \label{fig1}
\end{figure}
\end{center}

{\it Chebyshev polynomials of the second kind\/} are defined by
    $$U_r(\cos\theta)=\frac{\sin(r+1)\theta}{\sin\theta}$$
for $r\geq0$. Evidently, $U_r(x)$ is a polynomial of degree $r$ in
$x$ with integer coefficients. Chebyshev polynomials were invented
for the needs of approximation theory, but are also widely used in
various other branches of mathematics, including algebra,
combinatorics, number theory, and lattice paths (see \cite{K,
Ri}). For $k\geq0$ we define $R_k(x)$ by
$$R_k(x)=\frac{U_{k-1}\left(\frac{1}{2\sqrt{x}}\right)}
{\sqrt{x}U_k\left(\frac{1}{2\sqrt{x}}\right)}.$$ For example,
$R_0(x)=0$, $R_1(x)=1$, and $R_2(x)=1/(1-x)$. It is easy to see
that for any $k$, $R_k(x)$ is a rational function in $x$.

A {\em Dyck path} is a lattice path in the plane integer lattice
$\mathbb{Z}\times\mathbb{Z}$ consisting of up$\mn$steps $(1,1)$
and down$\mn$steps $(1,-1)$, which never passes below the $x$-axis
(see Figure~\ref{fig1}). Let $P$ be a Dyck path; we define the
{\em weight} of $P$ to be the product of the weights of all its
steps, where the weight of every step (up$\mn$step or
down$\mn$step) is $\sqrt{x}$. For example, Figure~\ref{fig1}
presents two Dyck paths, each of length $12$ and weight $x^6$.

A point on the Dyck path is called a {\em peak at height $k$} if
it is a point with coordinate $y=k$ that is immediately preceded
by a up$\mn$step and immediately followed by a down$\mn$step. For
example, Figure~\ref{fig1} presents two Dyck paths; the path on
the left has two peaks at height $2$ and two peaks at height $3$;
and the path on the right has one peak at height $1$, one peak at
height $2$, and one peak at height $3$. A point on the Dyck path
is called a {\em valley at height $k$} if it is a point with
coordinate $y=k$ that is immediately preceded by down$\mn$step and
immediately followed by up$\mn$step. For example, in
Figure~\ref{fig1}, the path on the left has two valleys at height
$1$ and one valley at height $2$, and the path on the right has
only two valleys at height $0$. The number of all Dyck paths
starting at $(0,0)$ and ending at $(2n,0)$ with exactly $r$ peaks
(resp. valleys) at height $k$ we denote by ${\rm peak}_k^r(n)$
(resp. ${\rm valley}_k^r(n)$). The corresponding generating
function is denoted by ${\rm Peak}_k^r(x)$ (resp.\ ${\rm
Valley}_k^r(x)$).

Deutsch~\cite{D} found the number of Dyck paths of length $2n$
starting and ending on the $x$-axis with no peaks at height $1$ is
given by the $n$th Fine number: $1$, $0$, $1$, $2$, $6$, $18$,
$57$, $\dots$ (see~\cite{D,DS,F} and \cite[Sequence M1624]{SP}).
Recently, Peart and Woan~\cite{PW} gave a complete answer for the
number of Dyck paths of length $2n$ starting and ending on the
$x$-axis with no peaks at height $k$. This result can be
formulated as follows.

\begin{theorem}{\rm(see \cite[Section~2]{PW})}\label{ppp}
The generating function for the number of Dyck paths of length
$2n$ starting and ending on the $x$-axis with no peaks at height
$k$ is given by
$$\frac{1}{1-\dfrac{x}{1-\dfrac{x}{1-\dfrac{\ddots}{1-\dfrac{x}{1-x^2C^2(x)} } } } },$$
where the continued fraction contains exactly $k$ levels.
\end{theorem}

Theorem~\ref{ppp} is in fact a simple consequence of
Theorem~\ref{rvv} (as we are going to show in Section 3).

\begin{theorem}{\rm(see \cite[Proposition~1]{RV})}\label{rvv}
For given a Dyck path $P$ we give every up$\mn$step the weight
$1$, every down$\mn$step from height $k$ to height $k-1$ not
following a peak the weight $\lambda_k$, and every down$\mn$step
following a peak of height $k$ the weight $\mu_k$. The weight of
$w(P)$ of the path $P$ is the product of the weights of its steps.
Then the generating function $\sum_P w(P)$, where the sum over all
the Dyck paths, is given by
$$
\frac{1}{1-(\mu_1-\lambda_1)-\dfrac{\lambda_1}{1-(\mu_2-\lambda_2)
-\dfrac{\lambda_2}{1-(\mu_3-\lambda_3)-\ddots}} }
$$
\end{theorem}

In this paper we find an explicit formulas for the generating
functions ${\rm Peak}_k^r(x)$ and ${\rm Valley}_k^r(x)$ for any
$k,r\geq0$. This allows us to express these functions via
Chebyshev polynomials of the second kind $U_k(x)$ and generating
function for the Catalan numbers $C(x)$. The main result of this
paper can be formulated as follows:

\begin{mtheorem}\
\begin{itemize}
\item[$(i)$] For all $k\geq2$,
          $${\rm Peak}_k^r(x)={\rm Valley}_{k-2}^r(x);$$

\item[$(ii)$] For all $k,r\geq0$,
          $${\rm Valley}_k^r(x)=\delta_{r,0}R_{k+1}(x)+\frac{x^rC^{r+1}(x)}
           {U_{k+1}^2\left(\frac{1}{2\sqrt{x}}\right)
           \Big(1-x(R_{k+1}(x)-1)C(x)\Big)^{r+1}};$$

\item[$(iii)$] For all $r\geq0$,
           $${\rm Peak}_1^r(x)=\delta_{r,0}+\frac{x^{3r+2}C^{2r+2}(x)}
            {(1-x^2C^2(x))^{r+1}}.$$
\end{itemize}
\end{mtheorem}
We give two proofs of this result. The first proof, given in
Section 2, uses a decomposition of the paths under consideration,
while the second proof, given in Section 3, uses the continued
fraction theorem due to Roblet and Viennot (see Theorem~\ref{rvv})
as the starting point.

\begin{remark}
By the first part and the second part of the Main Theorem, we
obtain an explicit expression for the generating function for the
number of Dyck paths starting at $(0,0)$ and ending on the
$x$-axis with no peaks at height $k\geq2$, namely
$${\rm Peak}_k^0(x)=R_{k-1}(x)+\frac{x^rC^{r+1}(x)}
           {U_{k-1}^2\left(\frac{1}{2\sqrt{x}}\right)
           \Big(1-x(R_{k-1}(x)-1)C(x)\Big)^{r+1}}.$$
\end{remark}

We also provide a combinatorial explanation for certain facts in
Main Theorem. For example, we provide a combinatorial proof for
the fact $(ii)$ in the Main Theorem for $r=k=0$.

\noindent{\bf Acknowledgments}. The author expresses his
appreciation to the referees for their careful reading of the
manuscript and helpful suggestions.
\section{Proofs: directly from definitions}
In this section we present a proof for the Main Theorem which is
based on the definitions of the Dyck paths.

\noindent{\bf Proof of the Main Theorem$(i)$}. We start by proving
the first part of the Main Theorem by introducing a bijection
$\Psi$ between the set of Dyck paths of length $2n$ with $r$ peaks
at height $k$ and the set of Dyck paths of length $2n$ with $r$
valleys at height $k-2$.

\begin{theorem}
${\rm Peak}_k^r(x)={\rm Valley}_{k-2}^r(x)$ for all $k\geq2$.
\end{theorem}
\begin{proof}
Let $P=P_1,P_2,\dots,P_{2n}$ be a Dyck path of length $2n$ with
exactly $r$ peaks at height $k\geq2$ where $P_j$ are the points of
the path $P$. For any point $P_j$ we define another point
$\Psi(P_j)=Q_j$ as follows. If $P_j$ appears as a point of a
valley at height $k-2$ then we define $Q_j=P_j+(0,2)$.  If $P_j$
appear as a point of a peak at height $k$ then we define
$Q_j=P_j-(0,2)$ (this is possible since $k\geq2$). Otherwise, we
define $Q_j=P_j$. Therefore, we obtain a new path
$Q=Q_1,Q_2,\dots,Q_{2n}$, and by definition of $Q$ it is easy to
see that $Q$ is a Dyck path of length $2n$ with exactly $r$
valleys at height $k-2$ (see Figure~\ref{bijpsi}).

\begin{center}
\begin{figure}[H]
\hspace*{2.55truecm} \epsfxsize=320.0pt \epsffile{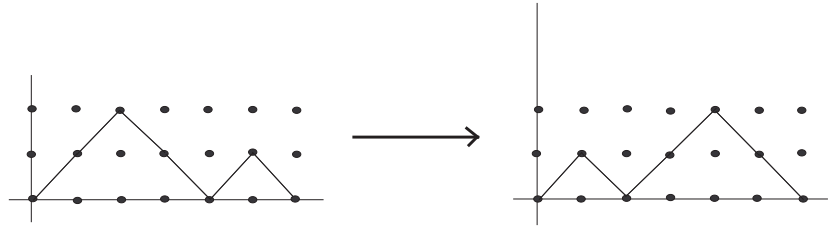}
\caption{Bijection $\Psi$.} \label{bijpsi}
\end{figure}
\end{center}

In fact, it is easily verified that the map which maps $P$ to $Q$
is a bijection. This establishes the theorem.
\end{proof}

\noindent{\bf Formula for ${\rm Valley}_0^0(x)$}. Let $P$ be a
Dyck path with no valleys at height $0$. It is easy to see that
$P$ has no valleys at height $0$ if and only if there exists a
Dyck path $P'$ of length $2n-2$ such that
                    $$P={\rm up\mn step}, P', {\rm down\mn step}.$$
Let $P''$ be the path that results by shifting $P'$ by $(-1,-1)$.
Then the map $\Theta$ which sends $P\rightarrow P''$ is a
bijection between the set of all Dyck paths starting at $(0,0)$
and ending at $(2n,0)$ with no valleys, and the set of all Dyck
paths starting at point $(0,0)$ and ending at $(2n-2,0)$. Hence
                     $${\rm Valley}_0^0(x)=1+xC(x),$$
where we count $1$ for the empty path, $x$ for the up$\mn$step and
the down$\mn$step, and $C(x)$ for all Dyck paths $P''$.

\noindent{\bf Proof of the Main Theorem$(ii)$}. First of all, let
us present two facts. The first fact concerns the generating
function for the number of Dyck paths from the southwest corner of
a rectangle to the northeast corner.

\begin{fact}\label{f11}{\rm (see \cite[Theorem A2 with Fact A3]{K})}
Let $k\geq0$. The generating function for the number of Dyck paths
which lie between the lines $y=k$ and $y=0$, starting at $(0,0)$
and ending at $(n,k)$ is given by
       $$F_k(x):=\frac{1}{\sqrt{x}U_{k+1}\left(\frac{1}{2\sqrt{x}}\right)}.$$
\end{fact}

The second fact concerns the generating function for the number of
Dyck paths starting at $(a,k+1)$ and ending at $(a+n,k+1)$ with no
valleys at height $k$.

\begin{fact}\label{f12}
The generating function for the number of Dyck paths starting at
$(a,k+1)$ and ending at $(a+n,k+1)$ with no valleys at height $k$
is given by
$$\frac{C(x)}{1-x(R_{k+1}(x)-1)C(x)}.$$
\end{fact}
\begin{proof}
Let $P$ be a Dyck path starting at $(k+1,0)$ and ending at
$(k+1,n)$ with no valleys at height $k$. It is easy to see that
$P$ has a unique decomposition of the form
  $$P=W_1,{\rm down\mn step},V_1,
  {\rm up\mn step},W_2,{\rm down\mn step},V_2,\dots, {\rm up\mn step},W_m,$$
where the following conditions holds for all $j$:
\begin{itemize}
\item[$(i)$] $W_j$ is a path consisting of up$\mn$steps and down$\mn$steps
starting and ending at height $k+1$ and never passes below the
height $k+1$;

\item[$(ii)$] $V_j$ is a path consisting of up$\mn$step and down$\mn$steps
starting and ending at height $k$ and never passes over the height
$k$ (see Figure~\ref{fig3}).
\end{itemize}

\begin{center}
\begin{figure}[H]
\hspace*{2.75truecm} \epsfxsize=325.0pt \epsffile{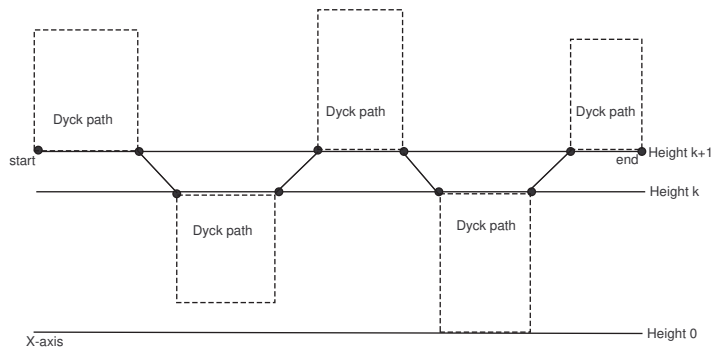}
\caption{A decomposition of a Dyck path starting at $(a,k+1)$ and
ending at $(a+n,k+1)$ with no valleys at height $k$.} \label{fig3}
\end{figure}
\end{center}

Using \cite[Theorem~2]{K} we get that the generating function for
the number of paths of type $V_j$ (shift for a Dyck path) is given
by $R_{k+1}(x)-1$. Using the fact that $W_j$ is a shift for a Dyck
paths starting and ending on the $x$-axis we obtain the generating
function for the number of Dyck paths of type $W_j$ is given by
$C(x)$. If we sum over all the possibilities of $m$ then we have
$$C(x)\sum_{m\geq0}\big(xC(x)(R_{k+1}(x)-1)\big)^m=\frac{C(x)}{1-x(R_{k+1}(x)-1)C(x)}.$$
\end{proof}

Now we are ready to prove the second part of the Main Theorem.

\begin{theorem}
The generating function ${\rm Valley}_k^r(x)$ is given by
$$\delta_{r,0}R_{k+1}(x)+\frac{x^{r}C^{r+1}(x)}
{U_{k+1}^2\left(\frac{1}{2\sqrt{x}}\right)\Big(1-x(R_{k+1}(x)-1)C(x)\Big)^{r+1}}.$$
\end{theorem}
\begin{proof}
Let $P$ be a Dyck path starting $(0,0)$ and ending at $(2n,0)$
with exactly $r$ valleys at height $k$. It is easy to see that $P$
has a unique decomposition of the form
$$\begin{array}{l}
P=E_1,{\rm up\mn step},D_0,{\rm down\mn step},{\rm up\mn
step},D_1, {\rm down\mn step},{\rm up\mn step},\dots,D_r,{\rm
down\mn step},E_2,
\end{array}$$ where the following conditions holds:
\begin{itemize}
\item[$(i)$] $E_1$ is a Dyck path that lies between the lines $y=k$ and
$y=0$, starting at $(0,0)$, and ending at point on height $k$;

\item[$(ii)$] $D_j$ is a Dyck path starting and ending at points on
height $k+1$ without valleys at height $k$, for all $j$;

\item[$(iii)$] $E_2$ is a Dyck path that lies between the lines $y=k$ and
$y=0$, starting at point on height $k$, and ending at $(2n,0)$.
\end{itemize}
Using Fact~\ref{f11} and Fact~\ref{f12} we get the the desired
result for all $r\geq1$. Now, if we assume that $r=0$, then we
must consider another possibility which is that all the Dyck paths
lie between the lines $y=k$ and $y=0$, starting at $(0,0)$, and
ending on the $x$-axis. Hence, using \cite[Theorem~2]{K} we get
that the generating function for the number of these paths is
given by $R_{k+1}(x)$.
\end{proof}

As a corollary of the Main Theorem$(ii)$ for $k=0$ (using
\cite[Example~1.18]{M}) we get

\begin{theorem}
For all $r\geq0$,
            $${\rm Valley}_0^r(x)=\delta_{r,0}+x^{r+1}C^{r+1}(x).$$
In other words, the number of Dyck paths starting at $(0,0)$ and
ending at $(2n,0)$ with exactly $r$ valleys at height $0$ is given
by
         $$\frac{r+1}{n}\binom{2n-r-1}{n+1}.$$
\end{theorem}

\noindent{\bf Proof of the Main Theorem$(iii)$}. If we merge the
first two parts of Main Theorem, then we get an explicit formula
for ${\rm Peak}_k^r(x)$ for all $r\geq0$ and $k\geq2$. Besides, by
definition there are no peaks at height $0$. Thus, it is left to
find ${\rm Peak}_1^r(x)$ for all $r\geq0$.

\begin{theorem}
For all $r\geq0$,
            $${\rm Peak}_1^r(x)=\delta_{r,0}+\frac{x^{3r+2}C^{2r+2}(x)}{(1-x^2C^2(x))^{r+1}}.$$
\end{theorem}
\begin{proof}
Let $P$ be a Dyck path starting at $(0,0)$ and ending at $(2n,0)$
with exactly $r$ peaks at height $1$. It is easy to see that $P$
has a unique decomposition of the form
$$P=D_0,{\rm up\mn step},{\rm down\mn step},D_1,{\rm up\mn step},
{\rm down\mn step}, \ldots,{\rm up\mn step}, {\rm down\mn
step},D_r,$$ where $D_j$ is a nonempty Dyck path starting and
ending at point on the $x$-axis with no peaks at height $1$.
Hence, the rest is easy to obtain by using \cite{D}.
\end{proof}

For example, for $r=0$ the above theorem yields the main result of
\cite{D}.
\section{Proofs: Directly from Theorem \ref{rvv}}
In this section we present another proof for the Main Theorem
which is based on Roblet and Viennot \cite[Proposition~1]{RV} (see
Theorem~\ref{rvv}).

Let $\lambda_j=x$ for all $j$, $\mu_j=x$ for all $j\neq k$, and
$\mu_k=z$. Theorem~\ref{rvv} yields
$$\sum_{r\geq 0}{\rm Peak}_k^r(x)z^r=\frac{1}{1-\dfrac{x}{1-\dfrac{x}{1-\dfrac{\ddots}{1- \dfrac{x}
{1-\dfrac{x}{1-(z-x)-\dfrac{x}{1-\dfrac{x}{1-\dfrac{x}{1-\ddots} }
} } } } } } },\eqno(1)$$ where $z$ appears in the $k$th level. On
the other hand, $xC^2(x)=C(x)-1$, we have that
$$C(x)=\frac{1}{1-\dfrac{x}{1-\dfrac{x}{1-\ddots}} }.\eqno(2)$$

Using the identities (1) and (2) with $xC^2(x)=C(x)-1$ we get

\begin{theorem}\label{aaa}
The generating function $\sum_{r\geq 0}{\rm Peak}_k^r(x)z^r$ is
given by
$$\frac{1}{1-\dfrac{x}{1-\dfrac{x}{1-\dfrac{\ddots}{1- \dfrac{x}
{1-\dfrac{x}{1-z-x^2C^2(x)} } } } } },$$ where the continued
fraction contains exactly $k$ levels.
\end{theorem}

For example, Theorem~\ref{aaa} yields for $z=0$ the generating
function ${\rm Peak}_k^0(x)$ as in the statement of
Theorem~\ref{ppp}. More generally, Theorem~\ref{aaa} yields an
explicit expression for ${\rm Peak}_k^r(x)$ for any $r\geq 1$ by
using the following lemma.

\begin{lemma}
For all $k\geq 1$,
$$\frac{1}{1-\dfrac{x}{1-\dfrac{x}{1-\dfrac{\ddots}{1- \dfrac{x}
{1-z-xA } } } }
}=R_k(x)\cdot\frac{1-zR_{k-1}(x)-xAR_{k-1}(x)}{1-zR_k(x)-xAR_k(x)}.$$
where the continued fraction contains exactly $k$ levels.
\end{lemma}
\begin{proof}
Immediately, by using the identity $R_{m+1}(x)=1/(1-xR_m(x))$ and
induction on $k$.
\end{proof}

Therefore, using Theorem~\ref{aaa}, the above lemma, and the
identity $R_{m+1}(x)=1/(1-xR_m(x))$, together with definitions of
$R_k(x)$, we get the explicit expression for the generating
function ${\rm Peak}_k^r(x)$ for any $r\geq1$ (see the Main
Theorem).


\begin{thebibliography}{WWW}
\bibitem[D]{D}
E.~Deutsch. Dyck path enumeration, {\em Disc. Math.} {\bf 204}
(1999) 167--202.

\bibitem[DS]{DS}
E.~Deutsch and L.W.~Shapiro, A survey of the Fine numbers, {\em
Disc. Math.} {\bf 241} (2001) 241--265.

\bibitem[F]{F}
T.~Fine, Extrapolation when very little is known about the source,
{\em Information and Control} {\bf 16} (1970), 331--359.

\bibitem[K]{K}
C.~Krattenthaler, Permutations with restricted patterns and Dyck
paths, {\em Adv. in Applied Math.} {\bf 27} (2001), 510--530.

\bibitem[M]{M}
S.~G.~Mohanty, Lattice path counting and applications, Academic
Press, 1979.

\bibitem[PW]{PW}
P.~Peart and W.J.~Woan, Dyck paths with no peaks at height $k$,
{\em J. of Integer Sequences} {\bf 4} (2001), Article 01.1.3.

\bibitem[RV]{RV}
E.~Roblet and X.G.~Viennot, Th\'eorie combinatoire des T-fractions
et approximants de Pad\'e en deux points, {\em Disc. Math.} {\bf
153} (1996) 271--288.

\bibitem[Ri]{Ri}
Th. Rivlin, Chebyshev polynomials. From approximation theory to
algebra and number theory, John Wiley, New York (1990).

\bibitem[SP]{SP}
N.J.A.~Sloane and S.~Plouffe, {\em The Encyclopedia of Integer
Sequences}, Academic Press, New York (1995).

\bibitem[S]{St}
R.~Stanley, Enumerative Combinatorics, vol. {\bf 1}, Wadsworth and
Brooks/Cole, Pacific Grove, CA, 1986, xi + 306 pages; second
printing, Cambridge University Press, Cambridge, 1996.
\end{thebibliography}
\end{document}